\documentclass{amsart}

\usepackage{epsfig}

\newcommand{\aaa}{{\mathcal A}}

\newcommand{\ddd}{{\mathcal D}}

\newcommand{\lll}{{\mathcal L}}
\newcommand{\ooo}{{\mathcal O}}
\newcommand{\eee}{{\mathcal E}}

\newcommand{\ppp}{{\mathcal P}}
\newcommand{\qqq}{{\mathcal Q}}

\newcommand{\rtrap}{{\mathcal{R}_{trap}}}
\newcommand{\resc}{{\mathcal{R}_{esc}}}

\newcommand{\xs}{{X_{S}}}

\newcommand{\dmax}{{d_{max}}}
\newcommand{\vmin}{{v_{min}}}

\newcommand{\R}{{\mathbb R}}
\newcommand{\Q}{{\mathbb Q}}

\newcommand{\Z}{{\mathbb Z}}
\newcommand{\C}{{\mathbb C}}
\newcommand{\N}{{\mathbb N}}

\newcommand{\vech}{{\vec{h}}}
\newcommand{\vecv}{{\vec{v}}}
\newcommand{\vecw}{{\vec{w}}}
\newcommand{\vecu}{{\vec{u}}}

\newtheorem{thm}{Theorem}[section]
\newtheorem*{thm*}{Theorem}
\newtheorem{cor}[thm]{Corollary}
\newtheorem{lem}[thm]{Lemma}

\newtheorem{dfn}[thm]{Definition}

\theoremstyle{definition}
\newtheorem{ex}{Example}

\begin{document}

\title[Substitutions and local finiteness]{Generalized $\beta$-expansions, substitution tilings, and local finiteness}
\author{Natalie Priebe Frank \and E. Arthur Robinson, Jr.}

\address{Natalie Priebe Frank\\Department of Mathematics\\Vassar College\\Box 248\\Poughkeepsie, NY  12604\\nafrank@vassar.edu}

\address{E. Arthur Robinson, Jr.\\Department of Mathematics\\George Washington University\\ Washington, DC 20052\\robinson@gwu.edu}

\begin{abstract}
For a fairly general class of
  two-dimensional tiling substitutions,  we prove that if the length
  expansion $\beta$ is a Pisot number, then the tilings defined
  by the substitution must be locally finite. We also give a simple example of a
  two-dimensional substitution on rectangular tiles, with a
  non-Pisot length expansion $\beta$, such that no tiling admitted by the
  substitution is locally finite. The proofs of both results
  are effectively one-dimensional and involve the idea of a certain
  type of generalized $\beta$-transformation.
\end{abstract}


\maketitle


\section{Introduction}
\begin{figure}[ht]
\epsfig{figure=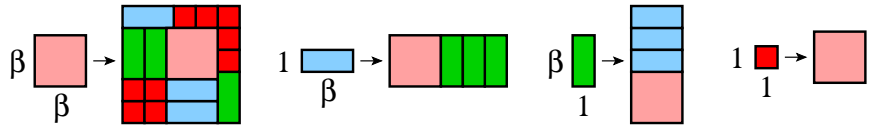,width=4.5in}
\caption{An inflate-and-subdivide rule $S$.}
\label{np.subs}
\end{figure}

\begin{figure}[ht]
\epsfig{figure=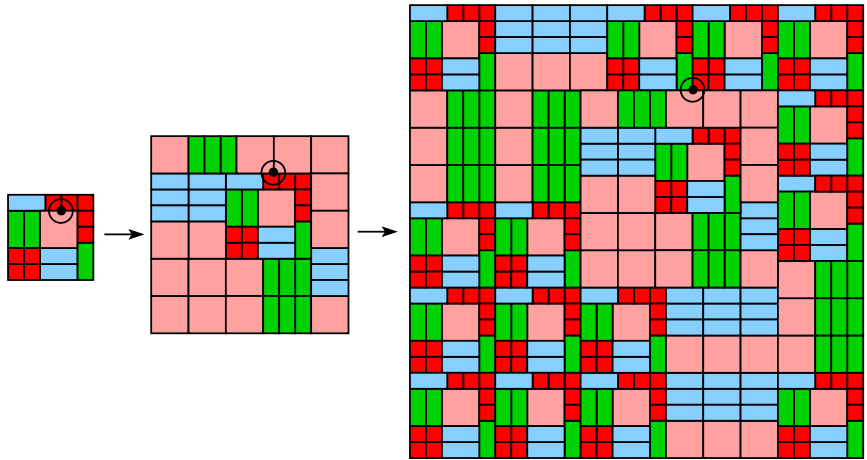,width=4.5in}
\caption{Three iterations of the larger square tile.}
\label{iterate.misfit}
\end{figure}

This paper was motivated by the study of simple tiling
substitutions such as the substitution $S$ shown  in Figure~\ref{np.subs}.
To the left of the arrows, four prototiles are shown along with their dimensions.  The larger
dimension, $\beta = (1+\sqrt{13})/2\approx 2.3027756377$, is a root of
the polynomial $p(x)=x^2 - x - 3$. The patches to the right of the
arrows show the application of the substitution $S$ to each prototile.
The substitution $S$ can be applied repeatedly, to produce arbitrarily
large patches of tiles. The first few applications, starting with the
large square prototile, are shown in Figure \ref{iterate.misfit}. 
Our interest in $S$ arises because its length expansion $\beta$ is not a
Pisot number.  An algebraic integer $\beta$ is called a
{\em Pisot number} if $\beta>1$, and all its Galois conjugates
$\beta'$ satisfy $|\beta'|<1$. 
Here, $\beta$ is not Pisot since its
 conjugate is $\beta' = (1-\sqrt{13})/2\approx -1.3027756377$.

It is known that whether or not the linear expansion of
 a tiling substitution $S$ is a Pisot number has implications for the
 tiling dynamical system that it defines.  (We make the distinction between the linear expansion and the length expansion of a substitution tiling precise in Section \ref{TilingResults}.)  Solomyak proved in
 \cite{Sol.self.similar} that a one- or two-dimensional self-similar tiling dynamical system with real linear expansion
 $\beta$ is weakly mixing if and only if $\beta$ is not Pisot.  A
 similar result holds for a two-dimensional self-similar
 tiling dynamical system with complex linear expansion
 $\lambda\in\C\backslash\R$. Solomyak showed that such a dynamical
 system is weakly mixing if and only if $\lambda$ is not a ``complex
 Pisot'' number.   However, in the two-dimensional
 case the proofs of these results, as well as many similar results, use
 the additional assumption of ``local finiteness'' (a property that is
 automatic for the one-dimensional case).  A tiling is defined to be
 {\em locally finite}\footnote{
This property is also sometimes also called
   ``finite local complexity''. In the theory Delone sets, however,
   ``local finiteness'' usually refers to a weaker property.}  
 if, up to rigid motions, it contains only finitely many two-tile
 configurations.  
 
 For substitution tilings, local finiteness follows from ``two-tile
 closure property''\cite{Robinson.ams}.  On the other hand, Danzer has
 conjectured \cite{danzer} that non-locally finite tilings may in fact
 be generic for certain classes of substitutions.  We will show that
 tilings admitted by the substitution in Figure~\ref{np.subs} fail to be
 locally finite.  It remains unclear the extent to which results like
 those in \cite{Sol.self.similar} can be extended to the non-locally
 finite case.

Although some examples of non-locally finite tilings have already
appeared (see \cite{danzer}, see also \cite{sadun}), our example has
the benefit of simplicity. It belongs to a family
\cite{frank} of tilings that can be created quite easily from any
one-dimensional substitution.  Moreover, our proof is based on a
one-dimensional method that may be of independent interest.

To see how  the lack of local finiteness appears from 
Figure~\ref{np.subs}, consider Figure~\ref{iterate.misfit}, which
shows three iterations of the large square tile.  In the first
iteration, we have circled a point where there is, in Danzer's
terminology \cite{danzer}, a ``misfit situation".
Under the inflate-and-subdivide rule $S$ the tiles are inflated by the
linear map $\beta I$; we can follow the circled point's location under
this linear map.  Following the circled point through two more
iterations, we see that the misfit it creates moves around relative to
the tile below it.  We will prove that this point occupies
infinitely many locations with respect to the tile below it.
Therefore the infinite tiling generated by this process is not
locally finite.

Now one can begin to see why the arguments are one-dimensional.
Inspection of Figure~\ref{np.subs} shows that there are two kinds of
upper edges: a long one and a short one, which we denote $l$ and $s$
respectively.  The reader can verify that a long upper edge is
substituted by a long followed by three short edges, ($l \to lsss$),
and a short upper edge is substituted by a long one, ($s \to l$).  The
circled point in Figure~\ref{iterate.misfit} is always on the corner
of two of the tiles above it, but it is (we will show) always in the
interior of an upper edge.  Thus we restrict our attention to the
one-dimensional substitution on the upper edges to prove our result.
In fact, following the misfit through these iterates is analogous to
applying a suitably defined ``generalized $\beta$-transformation",
which we discuss in the next section.

The organization of this paper is as follows.  We first focus on the
one-dimensional theory, which is itself of independent interest. In
Section \ref{Definitions},  we give the one-dimensional
definition of a substitution tiling and a generalized
$\beta$-transformation, and we show how they are related.  In Section
\ref{OneDResults} we give conditions that guarantee finite or infinite
orbits of the generalized $\beta$-transformations and relate this
result to one-dimensional tilings.  In Section \ref{TilingResults} we
apply this machinery to tilings of $\R^2$.

\section{Definitions} \label{Definitions}

\subsection{One-dimensional tiling substitutions}
 \label{def.one.d.subs}
 
 Let $\aaa=\{1, \dots, m\}$ and let $\aaa^*$ denote the set of all
 (nonempty) finite words in $\aaa$.  A {\em substitution} is mapping
 $\sigma:\aaa\rightarrow\aaa^*$.  Given a substitution $\sigma$, the
 {\em structure matrix} is the $d\times d$ matrix $M$ with entries
 $M_{i,j}$ equal to the number of occurrences of the symbol $j$ in the
 word $\sigma(i)$.  We require that $M$ has a real eigenvalue $\beta >
 1$ such that the right eigenvector $\vech=(h_1,h_2,\dots,h_m)$ has
 strictly positive entries.  It is worth noting that Perron-Frobenius
 theory ensures this in two special cases: irreducibility and
 primitivity.  We say $M$ is irreducible if for each $i, j$ there is
 an $n$ for which $(M^n)_{ij} > 0$.  If a single $n$ works for all
 $i,j$ then we say $M$ is primitive.  By the Perron-Frobenius theorem,
 if $M$ is irreducible then there exists an $r>0$ for which $\vech \in
 (r\Q(\beta))^m$.
  
 Given a substitution $\sigma$, we construct a one-dimensional tiling substitution $S$ and a corresponding tiling $\tau$
 of the line.  The {\em prototiles} $\ppp=\{P_1,\dots,P_m\}$ are
 intervals $P_j=[0,h_j)$ with lengths given by the entries of $\vech$.
A {\em tile} is a translation of a prototile
 $T=P_j+s=[s,s+h_j)$. A {\em tiling patch} $\pi$ is a finite sequence of
 tiles laid end to end:
\begin{equation}\label{patch}
\pi=\left(P_{j_1}+s_1,P_{j_2}+s_2,\dots,P_{j_n}+s_n\right)
\end{equation} 
where
\begin{equation*}
s_{i+1}=s_i + h_{j_i}.
\end{equation*}
Equivalently, in a tiling patch, the tiles are disjoint and the {\em
  support}, defined
\begin{equation}\label{support}
{\rm supp}(\pi):=\cup_{i=0}^n P_{j_i}+s_i,
\end{equation}
is an interval. A {\em tiling} $\tau$ is an infinite set of disjoint
tiles whose support is $\R$.

We call $s$ the left endpoint of the tiling patch $\pi$, and if ${\rm
  supp}(\pi)=[s,t)$, we call $t$ the right endpoint.  We can
specify a tiling patch $\pi$ 
by giving its left end point
$s$ and the word ${\bf w}=j_1 j_2 \dots j_n$ that specifies the prototiles
that occur in its sequence. In particular, the patch in
(\ref{patch}) is given by $\pi=\pi(s,{\bf w})$.

Now we can define the {\em tiling substitution} $S$ corresponding to
$\sigma$. It is the mapping from tilings to tiling patches defined by
\begin{equation*}
S(P_j+s)=\pi(\beta s,\sigma(j)).
\end{equation*} 
Clearly we have ${\rm supp}(S(T))=\beta\,{\rm supp}(T)$, and
$S(T+s)=S(T)+\beta s$. This allows us to extend $S$ to a mapping on
tiling patches, and even to tilings.  A tiling $\tau$ is said to be
$S$-admissible if every patch $\pi\subset\tau$ is equivalent to a
patch of $\pi'\subset S^n P$ for some $P\in\ppp$, $n\in\N$. The set of
all $S$-admissible tilings (which is always nonempty) is called a {\em
  tiling space}, and is written $X_S$.

\begin{ex} \label{OneDSubs}
  Consider the substitution $ l \to lsss ,\, s \to l$ from the
  introduction. The structure matrix is
  $$M = \left[ \begin{array}{cc} 1&3\\1&0 \end{array}\right],$$
  which has characteristic polynomial $p(x) = x^2 - x - 3$, and the
  Perron eigenvalue $ \beta = (1 + \sqrt{13})/2$. Taking $\vech =
  (\beta, 1)$ as a positive right eigenvector, the prototiles are $P_1
  = [0, \beta)$ and $P_2 = [0, 1)$.  These tiles have the same lengths
  as the edges of the tiles in Figure~\ref{np.subs}.  We will continue
  our discussion of this example below.
\end{ex}




\subsection{Generalized $\beta$-transformations}
\label{Gen.beta}
Fix an interval $E = [0, B)$
and a real number $\beta>1$.

\begin{dfn}
  We say a mapping $F_\beta : E \to E$ is a {\em generalized
    $\beta$-transformation} if (i) $F_\beta$ is piecewise linear,
  (ii) $F_\beta$ is continuous from the right and (iii)
  $F'_\beta(x)=\beta$ for almost every $x\in I$. 
\end{dfn}

More explicitly, $F_\beta$ is a generalized $\beta$-transformation if
there is a partition of $E$ into $l$ subintervals $I_j=[x_{j-1},x_j)$,
with endpoints $0 = x_0 < x_1 < ... < x_l=B$, and real numbers $y_1,
y_2, ...  y_l$ with $\beta I_j - y_j \subseteq E$ for $j=1,\dots,l$, so that 
\begin{equation}\label{beta.def}
F_\beta(x) = \beta x - y_j \text{ whenever }  x \in I_j. 
\end{equation}

One example of this is the classical $\beta$-transformation
$T_\beta:[0,1)\rightarrow[0,1)$, defined $T_\beta(x)=\beta x\text{ mod
}1$ (see Renyi~\cite{renyi}). A type of generalized
$\beta$-transformation, more general than what we consider here, was
studied by Wilkinson~\cite{wilkinson}.
He showed under several
additional conditions that $F_\beta$ has an absolutely continuous
invariant measure and is weakly Bernoulli with respect to the
partition into the intervals $I_j$.

\begin{dfn}
We call a generalized $\beta$-transformation $F_\beta$ {\em algebraic} if 
$\beta$ is an algebraic number and there is an $r \in \R$ for which 
\begin{equation*}
x_0,\dots,x_l,y_1,\dots,y_l\in r\Q(\beta).
\end{equation*}
\end{dfn}

\subsection{Substitution tilings and generalized $\beta$-transformations}
\label{subs.and.beta}
\begin{figure}[ht]
\epsfig{figure=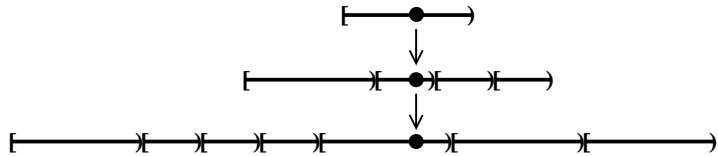,height=.8in}
\caption{The map $\bar S$ finds the central 
tile after the substitution is applied.}
\label{one.d.iterate}
\end{figure}

Now let us think of the circled point in Figure \ref{iterate.misfit}
as being the origin situated in the interior of the
one-dimensional tile that is the upper edge of a two-dimensional tile
(Figure~\ref{one.d.iterate}). 

Suppose that $T=[s,t)$ is a tile with $s\le 0$ and $t>0$, so that
$0\in T$. We call such a tile a {\em central tile}, illustrated in
Figure~\ref{one.d.iterate} with $T = [1-\beta, 1)$. It follows that
$0\in {\rm supp}(S(T))=[\beta s,\beta t)$, and thus there is a unique
central tile in
$S(T)=\left(P_{j_1}+s_1,P_{j_2}+s_2,\dots,P_{j_n}+s_n\right)$, i.e.,
$0\in P_{j_i}+s_i$ for some $i$. We can therefore define a
mapping $\bar S$ from the set of all central tiles to itself by taking
$\bar S(T)$ to be the central tile in $S(T)$. We are
interested in the effects of iterating $\bar S$ on central tiles, but it
is convenient to think of this process a little differently.

First we consider the patch $\pi:=\pi(0,123\dots
m)=(E_1,E_2,\dots,E_m)$ where $E_j=J_j+g_j$ with $g_0=0$ and
$g_{j+1}=g_j+h_j$. We write the support of $\pi$ as $E=[0,B)$.
The mapping $F : E \rightarrow E$ is defined as follows.  For $x \in E$,
find $j \in \aaa$ so that $x \in E_j$. Then $E_j - x$ is a central 
tile, and 
${\bar S}(E_j-x)=E_k-y$ for some
unique $k \in \aaa$ and $y \in \R$. We put $y=F(x)$, noting that $F(x)
\in E_k$.

We may write $F$ explicitly as a generalized $\beta$-transformation by
partitioning each $E_j$ into subintervals that will be sent by $F$
onto whole tiles (the canonical subdivision suggested by $S$).  We
obtain a partition $0=x_0 < x_1 < ... < x_l = B$ into intervals $I_1,
..., I_l$ such that for each $I_k$ there is a $j$ such that $F(I_k) =
E_j$.  The restriction of $F$ to $I_k$ is the linear map taking $I_k$
onto $E_j$. So for $x \in I_k$ we have
\begin{equation}
F(x) = \beta(x - x_{k-1}) + g_{j-1};
\end{equation}
hence $y_k = -\beta x_{k-1} + g_{j-1}$.

\begin{ex} \label{BetaEx}
The generalized $\beta$-transformation for Example~\ref{OneDSubs} is defined on the interval $E
  = [0,1+\beta)$ so that $E_1 = [0, \beta)$ and $E_2 =
  [\beta, \beta + 1)$.  The partition of $E$ given by the
  substitution is $x_0 = 0, x_1 = 1, x_2 = 1+1/{\beta}, x_3 = 1 +
  2/{\beta}, x_4 = 1 + 3/{\beta} = \beta,$ and $x_5 = 1+\beta$.  Since
  $g_0 = 0$ and $g_1 = \beta$ our translations are $y_1 = 0, y_2 =
  0,y_3 = 1, y_4 = 2,$ and $y_5 = 3 + \beta$.  The transformation is
  as shown in Figure \ref{beta.ex}.
\begin{figure}[ht]
\epsfig{figure=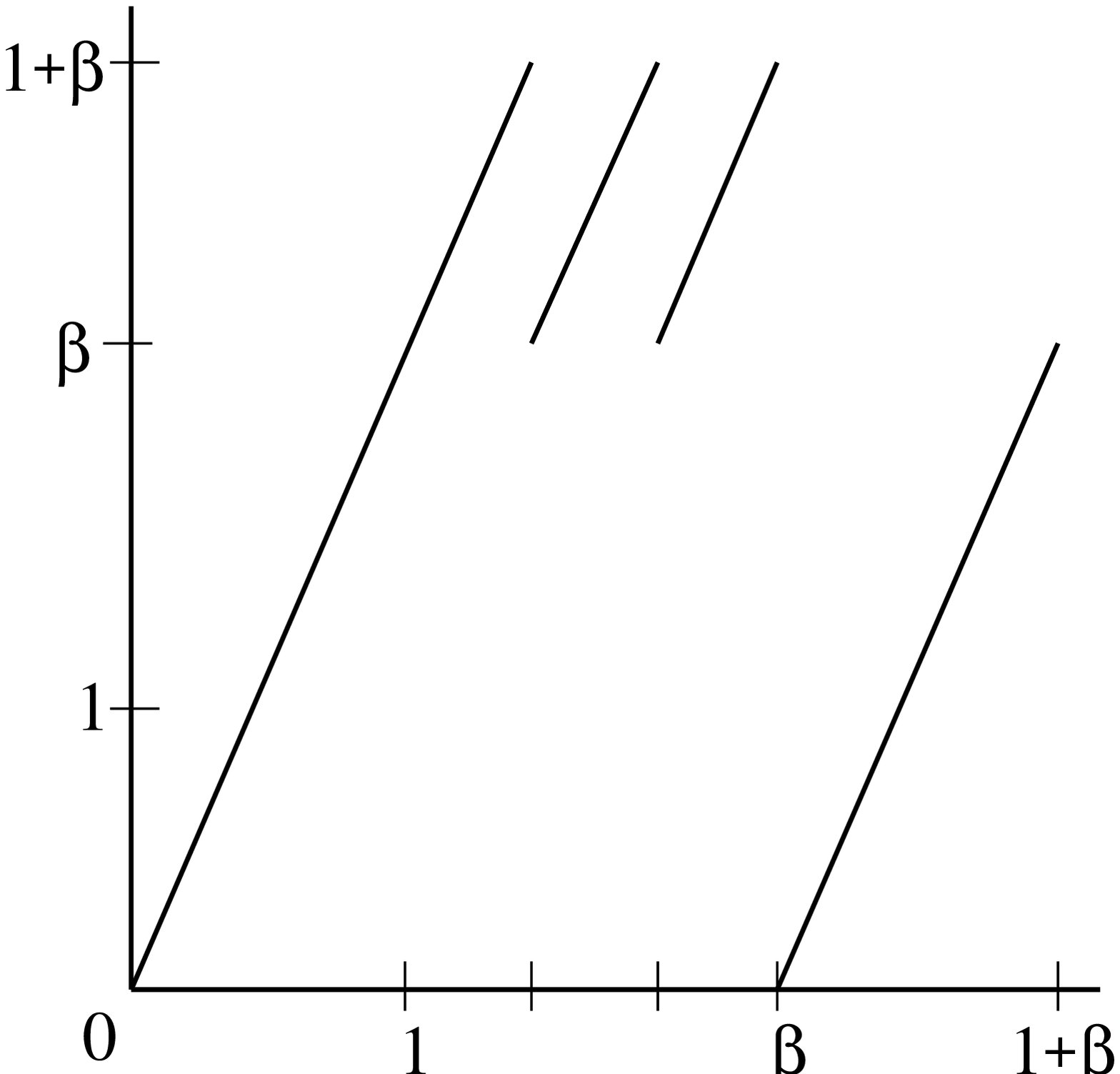,width=2.5in}
\caption{The map $F: E \to E$.}
\label{beta.ex}
\end{figure}

\end{ex}

Since each interval in a tiling substitution expands by a factor of $\beta$, we have the following lemma.

\begin{lem}
  The mapping $F$ given by a tiling substitution with expansion
  constant $\beta$ is an algebraic generalized $\beta$-transformation whenever $\vec{h} \in (r\Q(\beta))^m$.
\end{lem}

This is certainly true whenever $M$ is irreducible or primitive.


\section{One-dimensional results}
\label{OneDResults}


Suppose $\beta > 1$ is a root of the irreducible monic polynomial
$p(x)= g_0 + g_1 x + g_2 x^2 + ... + x^d\in \Z[x]$, where $d>1$.  We
can identify $\Q(\beta)$ with $\Q^d$ via the basis $\{1, \beta,
\beta^2, ... \beta^{d-1}\}$ over $\Q$.  Multiplication by $\beta$ in
$\Q(\beta)$ corresponds to multiplication by $A$ in $\Q^d$, where
\begin{equation} A =
\left[\begin{array} {cccc}0 & 0 & \cdots & -g_0\\1 & 0 & \cdots & -g_1\\0 & 1 & \cdots & -g_2\\
\cdots&&&\\0&  \cdots & 1 & -g_{d-1}\\
\end{array}\right]
\end{equation}
is the companion matrix of
$p(x)$.

To fully understand the action of $A$, we embed $\Q^d$ in $\C^d$
(although if the eigenvalues are all real, it suffices to embed in
$\R^d$).  
We define a surjective linear map $\phi:\C^d\rightarrow \C$ by  
$$\phi(\vecw) = w_0 + w_1 \beta + w_2 \beta^2 + ... + w_{d-1}
\beta^{d-1}.$$
The restriction of $\phi$ to $\Q^d$ is the isomorphism
between $\Q^d$ and $\Q(\beta)$.  Moreover,  for any $\vecw \in
\C^d$,
\begin{equation} \label{mult.by.A}
\phi(A \vecw) = \beta \phi(\vecw).
\end{equation}

The Galois conjugates $\beta_2, \beta_3, ... \beta_d$ of $\beta$ are
eigenvalues of $A$ with eigenvectors $\vecv_2, \vecv_3, ...  \vecv_d$.
By linearity $\phi(A \vecv_j) = \phi( \beta_j \vecv_j) = \beta_j
\phi(\vecv_j)$, but by (\ref{mult.by.A}) we also have $\phi( A
\vecv_j) = \beta \phi(\vecv_j)$, and so $\phi(\vecv_j) = 0$ for all $j
= 2, 3, ..., d$.  Let $V_0$ be the subspace of $\C^d$ spanned by
$\vecv_2,\dots,\vecv_d$, and let $V^+$ be the (one-dimensional)
eigenspace of $\beta$.  Since $\phi$ is not the zero map, we can
select the vector $\vecv \in V^+$ with $\phi(\vecv) =1$ as a basis
vector.

Consider fixed  an algebraic
generalized $\beta$-transformation $F : E \to E$.  Recalling the
notation of Subsection \ref{Gen.beta} we have the partition of $E$
given by $0 = x_0 < x_1 < ... < x_l=B$ with $x_j \in \Q(\beta)$; we
note that by definition $\beta (x_j- x_{j-1}) \le B$.  We denote the
$j$th subinterval as $I_j = [x_{j-1}, x_j)$, and we have
$y_1, y_2, ... y_l \in \Q(\beta)$ so that $\beta I_j - y_j \subseteq
E$.

Since $\phi(\vecu) = 0$ for any $\vecu \in V_0$, we know that $\phi(c
\vecv + \vecu) \in \R$ if and only if $c \in \R$.  Thus we can define
$$
W_{[a,b)} = \{\vecw : \phi(\vecw) \in [a,b) \} = \{c \vecv + \vecu:
\vecu \in V_0 \text{ and } c \in \left[a, b \right)\}\subseteq\C^d
$$
to be the set of all points in $\C^d$ that $\phi$ maps onto the
real interval $[a,b)$.

Let $W_E$ denote the subset of $\C^d$ corresponding to $E$ and
denote the subsets of $\C^d$ that correspond to the intervals $I_j$ as
$W_j = W_{I_j}$; the $W_j$s form a partition of $W_E$.  There are unique vectors $\vecw_1, \vecw_2, ...,
\vecw_d$ in $\Q^d$ for which $\phi(\vecw_i) = y_i$, and one can easily
check that $A W_j - \vecw_j \subseteq W_E$.
Thus we can define the map $f: W_E \to W_E$ as
\begin{equation} 
 f(\vecw) = A \vecw - \vecw_j \text{ whenever }  \vecw \in W_j \label{fff.eqn};
\end{equation}
and we see that
\begin{equation} \phi \circ f = F \circ  \phi. \end{equation}

Because our interest lies with the mapping $F$, the theorem below is
stated in terms of it, but the result is actually proved for $f$ and
then extended to $F$.  The proof is mostly linear algebra and uses
ideas similar to \cite{danzer}, \cite{schmidt} and \cite{thurston} .
\begin{thm} \label{orbit.type}
Let $\beta, \beta_2, ..., \beta_d$, $E$, and $F$ be defined as above.
\begin{enumerate}
\item If $|\beta_j| < 1$ for all $j = 2, 3, ..., d$, then for any $x
  \in \Q(\beta)\cap E$ the orbit of $x$ under $F$ is finite and
  therefore eventually periodic. \label{pisot.finite}
\item If $|\beta_j| > 1$ for some $j \in \{2, 3, ..., d\}$, then there
  exists an $x \in \Q(\beta)\cap E$ such that the orbit of $x$ under
  $F$ is infinite. \label{nonpisot.infinite}
\end{enumerate}
\end{thm}

\begin{proof}
  Since $\phi \circ f = F \circ \phi$ and the map $\phi$ is a
  bijection when restricted to $\Q^d$ it suffices to prove the results
  for $f$.
  
  Suppose $|\beta_j| < 1$ for all $j = 2, 3, ..., d$. Let $\vecw \in \Q^d \cap W_E$ and $\vecw_i$ from
  (\ref{fff.eqn}) be written in terms of the eigenvectors of $A$ as
  $$\vecw = c \vecv + \sum_{i = 2}^d c_i \vecv_i \quad \text{ and }
  \quad \vecw_j = w_{j,1} \vecv + \sum_{i = 2}^d w_{j,i} \vecv_i. $$
For $i=2,\dots,d$ choose
$$
  r_i > \frac{|w_{j,i}|}{1 - |\beta_i|} \quad \text{ for all }
  \quad j \in 1, 2, ..., d.
  $$
  Define the {\em trapping region} $\rtrap$ to be the set of all
  points in $W_E$ such that $c_i \le r_i$.  

We claim that $\rtrap$ is
  invariant under the action of $f$.  Suppose $\vecw \in \rtrap$ with
  such a choice of values, and find $j \in 1, 2, ..., l$ such that  $\vecw \in W_j$.  Then
\begin{equation}
f(\vecw) = A \vecw - \vecw_j = 
(\beta \vecv + \sum_{i = 2}^d \beta_i c_i \vecv_i) -
(w_{j,1} \vecv + \sum_{i = 2}^d w_{j,i} \vecv_i).\label{fff(vecw)}
\end{equation}
Since the component in the direction of $\vecv$ is unrelated to the
definition of $\rtrap$ we see that we must have $|\beta_i c_i -
w_{j,i}| \le r_i$ for $i = 2, 3, ... d$.  This is true since
$|w_{j,i}| < (1 - |\beta_i|)r_i $ and $|c_i| \le r_i$.  This proves 
the invariance. 

Now note that there is some $q \in \Z$ for which $\vecw$ and the
$\vecw_j$s are in $\frac{1}{q} \Z^d$.  Since $A$ is an integer matrix
it leaves $\frac{1}{q} \Z^d$ invariant and so $f\left(\frac{1}{q} \Z^d
  \cap W_E\right) \subseteq \frac{1}{q} \Z^d \cap W_E$.  Choose the
numbers $r_i$ for the trapping region so that $\vecw \in \rtrap$; clearly
$f\left(\frac{1}{q} \Z^d \cap \rtrap\right) \subseteq \frac{1}{q} \Z^d
\cap \rtrap$.  Since $\frac{1}{q} \Z^d \cap \rtrap$ is a finite set
the orbit of $\vecw$ is finite and thus eventually periodic.

To prove the second part of the theorem we must show there is an
infinite {\em escape region} $\resc \subset W_E$ for which orbits
beginning in $\resc$ tend to infinity.  Without loss of generality we
may assume that $|\beta_2| > 1$.  Choose $r_2$ so that
$$ r_2 > \frac{|w_{j, 2}|}{|\beta_2| - 1} \quad \text{ for all } \quad j = 1, 2, ..., d.$$
We define
$$\resc = \{ \vecw = c \vecv + \sum_{i = 2}^d c_i \vecv_i \subset W_E
\text{ such that } |c_2| > r_2\}.$$

To prove the result, we must show that for any $\vecw \in \resc \cap
\frac{1}{q} \Z^d$, the component of $f(\vecw)$ in the direction of
$\vecv_2$ exceeds that of $\vecw$ in modulus.  Supposing again that
$\vecw \in W_j$, by the computation in equation (\ref{fff(vecw)}) this
amounts to showing that $|\beta_2 c_2 - w_{j,2}| > |c_2|$.  By
choice of $r_2$ we have that $|w_{j,2}| < (|\beta_2| - 1) r_2 \le
(|\beta_2| - 1) |c_2|$ and the result follows.

\end{proof}

\begin{cor}
  Let $S$ be a one-dimensional tiling substitution with expansion constant $\beta>1$ and $\vech \in (r\Q(\beta))^m$.  Let $T$ be a central tile.
  (i) Suppose $\beta$ is a Pisot number.  If the endpoints of $T$
  belong to $\Q(\beta)$, then there exist $K,L\in \N$ so that ${\bar
    S}^{NL+K}(T)={\bar S}^{K}(T)$ for all $N>0$. (ii) Suppose $\beta$ is
  non-Pisot.
  Then there is a $q\in \N$ and finite 
  $E\subseteq\frac{1}{q}\Z[\beta]$ such that for any
  $s\in\frac{1}{q}\Z[\beta]\backslash E$, if $T=[-s,t)$ is central, then
  ${\bar S}^{K}(T)\not={\bar S}^{L}(T)$ for any $K>L\ge 1$.

\end{cor}

\section{Tiling Results}
\label{TilingResults}

In this section we study a class of self-similar substitution tilings
in $\R^2\cong\C$. Our first result is that if a substitution $S$ has a
Pisot length expansion, then the tilings it produces are
locally finite.  This solves a conjecture of Danzer \cite{danzer} for the type of tilings that we consider.  
In this section we also prove that
tilings generated by the substitution from Figure~\ref{np.subs} are
not locally finite.

We begin with some 
general definitions.
A {\em tile} $T$ is a closed topological disk $D\subseteq\C$, together
(possibly) with a label which is an element of some fixed finite set
$\aaa$. We say $D={\rm supp}(T)$ is the {\em support} of $T$. Two
tiles are called {\em equivalent} if they carry the same label and
have congruent supports.

We will  fix a finite set $\ppp$ of inequivalent tiles, called
{\em prototiles}, and assume all tiles $T$ are equivalent
to tiles from $\ppp$. A {\em packing} is a collection tiles
that intersect only on their boundaries, and the support of a packing
is the union of the supports of its tiles\footnote{
In the one dimensional case, it was
  convenient to take tiles as half-open and insist on disjointness.
  This also forced equivalence to preserve orientation, a restriction
  we do not impose here.}.
A {\em tiling patch} $\pi$ is a finite packing with a connected
support. We denote the set of patches by $\ppp^*$. A {\em tiling}
$\tau$ is a packing with support $\C$.

Let $\lambda \in \C$
with $|\lambda| > 1$.  Suppose $S: \ppp \to \ppp^*$ is such that for
each $P \in \ppp$, the ${\rm supp}(S(P))$ is congruent\footnote{
In the
  event that there is more than one congruence between the supports
  $S(P)$ and $\lambda\cdot P$, we will use the implicit
  parameterization of the boundaries to ensure uniqueness of $S(T)$.
} 
to $\lambda\cdot{\rm supp}(P):=\lambda\cdot P$. We call $S$ a {\em
  self-similar tiling substitution} (or simply a {\em substitution}
when there is no danger of confusion). We call $\lambda$ the {\em linear 
expansion} of $S$, equating multiplication by $\lambda$ with a linear
map on $\R^2$; volumes must expand by a factor of $\lambda \overline{\lambda}$.   We call $\beta=|\lambda|$ the {\em length expansion} of $S$ because it gives the expansion of one-dimensional objects. 

 The substitution $S$ extends in an obvious way to any
tile $T$, and can thus be iterated to obtain infinitely large patches
$S^n(P)$.  Following the one-dimensional case, we say a tiling $\tau$ is
{\em $S$-admissible} if for every patch $\pi\in\tau$ is equivalent to
a patch $\pi'\in S^n(P)$ for some $P\in\ppp$ and $n\ge 0$. We denote
the set of $S$-admissible tilings $\tau$ by $X_S$.

The following Lemma may be well-known, but we include a short proof for
completeness.
\begin{lem} \label{sst.no.circles}
  Let $\ppp$ be a finite prototile set and let $S$ be a tiling
  substitution on $\ppp$.  No prototile in $\ppp$ contains a nontrivial circular arc
  in its boundary. \end{lem}

\begin{proof}
  Suppose the boundary of a prototile $P \in \ppp$ contains an arc of
  radius $r>0$ subtending an angle $\theta >0$.  The boundary of
  $S^n(P)$ contains arcs of radius $\beta^nr$ subtending the angle
  $\theta$, and this implies that there are prototiles in $\ppp$ that
  have circular arcs of radius $\beta^nr$ in their boundaries.
  Since the prototile set is finite, we can find the maximum angle
  $\theta_n>0$ subtended by the arc of radius $\beta^nr$ in any
  prototile.  The arclengths $ \theta_n\beta^nr$ must go to zero
  in order for the prototiles to be topological disks.
  
  Let $\dmax$ and $\vmin$ be the maximum diameter and minimum volume
  of prototiles in $\ppp$, respectively.  The tiles inside of $S^n(P)$
  adjacent to the circular arc of length $\theta \beta^n r$ lie
  inside a strip of volume
  $$\theta/2[(\beta^n r)^2 - (\beta^n r - \dmax)^2] +
  2\pi\dmax^2,$$
  where $2\pi\dmax^2$ ensures the ends of the strip are
  covered.  The maximum number of tiles that can fit inside this strip
  is
$$\frac{\theta/2[2\dmax \beta^n r -\dmax^2] + 2\pi\dmax^2}{\vmin}
< \frac{\theta \dmax \beta^n r }{\vmin}+K,$$
where $K=(2\pi -
\theta/2)\dmax^2/\vmin$ does not depend on $n$.  The minimum number of
tiles required to cover the arc of length $\theta \beta^n r $ is
$(\theta\beta^n r )/(\theta_n \beta^n r )=\theta/\theta_n.$
For all $n = 1, 2, ...$ we must have
$$\frac{\theta}{\theta_n}\le \frac{\theta \dmax \beta^n r }{\vmin}
+K , \text{ or }$$
$$\theta \le \frac{\theta \dmax}{\vmin} (\theta_n \beta^n r) + K
\theta_n.$$
We have noted that both terms on the right must go to
zero, forcing $\theta$ to be zero also.  This contradiction completes
the proof.
\end{proof}

\subsection{Local finiteness when the length expansion is Pisot}
\label{PV.loc.finite}

For the result we are going to prove, we will need some mild
additional geometric hypotheses on prototiles and tilings.  A {\em
  straight edge} $\ell$ in the boundary of a tile $T$ is a line
segment that cannot be extended without leaving the boundary.  Any
such edge has a left and a right {\em endpoint}, relative to the
boundary being traversed clockwise.  The absence of circular arcs
explains the terminology in the following definition.
\begin{dfn}
A tile $T$ is called {\em fractagonal} if: 
\begin{enumerate}
\item $T$ has only finitely many straight edges, and
\item for any straight edge $\ell$ of $T$, the rest of $T$ lies
  in the interior of a half-plane generated by $\ell$ (we call this
  {\em weak convexity}).
\end{enumerate}
\end{dfn}
\noindent The authors conjecture that the first property is
  automatically 
satisfied by the prototiles of a self-similar tiling substitution, provided they are weakly convex.

We also need an assumption that for straight edges, collinearity
implies comeasurability in the following sense:
\begin{dfn}  
  We say a tiling substitution $S$ is a {\em property-(C) substitution} if
  whenever two straight edges $\ell_1$ and $\ell_2$ are collinear in a
  tiling $\tau \in X_S$, there is an $r>0$ for which $|\ell_1|,
  |\ell_2| \in r \Q(\beta)$, where $\beta=|\lambda|$ is the length expansion of $S$.
\end{dfn}

\begin{thm} \label{thm.flc}
  Suppose that $S$ is a property-(C) tiling substitution on a finite
  set $\ppp$ of fractagonal prototiles.
If the length expansion of $S$
  is a Pisot number, then any tiling in $\tau \in \xs$ is locally
  finite.
\end{thm}

\noindent{\bf Notes.}   (1) We point out two common
assumptions that are not necessary here.  The substitution need not be
primitive in that the structure matrix $M$ may never satisfy $M^n >
0$.  Tilings in $\xs$ need not be translationally finite, i.e.
prototiles may appear in an infinite number of orientations.  

(2) The authors are not aware of any primitive tiling substitutions that are
not property-(C).

(3) In his dissertation \cite{frettloeh}, Dirk Frettl\"oh obtained a similar result under much stronger assumptions, including that the substitution is primitive, the expansion is an integer, and that the tiles are polygonal.

\begin{proof}
  A {\em simple adjacency} in $\tau \in \xs$ is a two-tile patch
  $\alpha\subset \tau$ where the intersection is along straight edges.
  Kenyon \cite{kenyon} proves that for any tiling $\tau$ made from a
  finite prototile set $\ppp$, if $\tau$ has an infinite number of
  inequivalent two-tile patches, then those patches are either simple
  adjacencies, or they occur along an entire circle of tile
  boundaries.  By Lemma~\ref{sst.no.circles} the latter cannot occur
  in a substitution tiling.  It follows that in order to show
  $\tau\in\xs$ is locally finite, it suffices to show that there are
  only finitely many simple adjacencies $\alpha$.
  
  
  Fix a simple adjacency $\alpha$ in $\tau \in \xs$.  Two straight
  edges define $\alpha$; the difference between two vertices, one
  selected from each edge, is called an {\em offset} of $\alpha$.
  Thus each $\alpha$ defines four different offsets. It suffices to
  show that the total number of offsets is finite.
  
  Let the offset $\ooo$ of $\alpha$ be fixed.  We want to see $\ooo$
  as the difference of sums of edge lengths of tiles (see
  Figure~\ref{sum.powers}).  As such, it is convenient to consider
  edges as being intervals on the positive real axis.  We call a tile
  $T$ placed below (resp. above) the real axis with a straight edge on
  the real axis a {\em lower} (resp.  {\em upper}) tile.  By applying
  a rigid motion, we can take edge representatives with their left
  (resp. right) endpoints at the origin. We call these {\em edge sets}
  $\eee_\ell$ and $\eee_u$.  
  
By weak convexity, the substitution $S$ induces a one-dimensional substitution $S_\ell$ on $\eee_\ell$ in the sense of Subsection~\ref{def.one.d.subs}.    It is clear that the length expansion $\beta$ of
$S$ is the expansion of $S_\ell$, so the one-dimensional substitution has a Pisot expansion constant.
%
   If $T$ is the lower tile containing $e \in \eee_\ell$, we simply
  rotate $S(T)$ so that $\lambda\cdot e$ is mapped to $\beta \cdot e$
  and look at the lower tiles of $\beta\cdot e$.  Similarly we can
  create the one-dimensional substitution $S_u$ on $\eee_u$ by
  considering the substitution on upper tiles.  In particular, $S_u$
  is just the reverse of $S_\ell$.
  
 By property (C), the edge sets $\eee_\ell$ and $\eee_u$ are partitioned into a finite
  number of subsets for which there are constants $r_1, ..., r_k > 0$ with each $r_i/r_j$, $i\not=j$, irrational modulo $\beta$, such that
  each edge length is in $r_i \Q(\beta)$ for some $i$.    Note that the application of $S_\ell$ or $S_u$ preserves these classes.   Without loss of generality we may assume that the edges of the simple adjacency $\alpha$ have lengths in $\Q(\beta)$ (i.e. $r = 1$).
  
  By the construction of $\tau$, there is some smallest $n \ge 0$ for
  which there exists a tile $T$ such that the simple adjacency $\alpha$ is a sub-patch of
  $S^{n+1}(T)$.  The intersection of the two tiles in $\alpha$ is
  contained in a line segment $L$ that is composed of tile edges and
  is maximal in the sense that it cannot be extended without either
  entering the interior of some tile in $S^{n+1}(T)$ or leaving $S^{n+1}(T)$
  entirely.  We call $L$ a {\em fault line} through the tile $S^{n+1}(T)$.
  Since $\ooo$ depends only on $\alpha$, we can use $L$ to compute the
  value of $\ooo$, even if the tiles surrounding $L$ in $S^{n+1}(T)$ do
  not appear in $\tau$ (except of course those in $\alpha$).

  If $n=0$, then $S(T)$ can be moved to the positive $x$-axis so
  that an endpoint of $L$ lies at the origin.  By the weak convexity
  assumption it follows that $L$ can be seen as a union of edges of
  upper tiles and also as a union of edges of lower tiles.  In this
  case $L$, along with the  information about which level-one tile
  it is in and which tile edges define it, is called an {\em initial
    segment}.
  
  If $n>0$, then by the minimality of $n$, there must be two
  level-$n$ tiles whose boundaries intersect $\alpha$.  From this
  we see that $L$ must contain a line segment $L'$ that is a
  union of edges of level-$n$ tiles.  Since the tiles are assumed
  to be fractagonal, weak convexity implies that if $S^{n+1}(T)$ is moved
  so that $L'$ in on the $x$-axis, then $L'$ is the union of edges of
  upper level-$n$ tiles and also the union of edges of lower
  level-$n$ tiles.  Thus it is the image under $S^{n}$ of some
  initial segment $L''$.  
    
  Since there are only a finite number of level-one tiles up to
  equivalence, there are only a finite number of initial segments up
  to equivalence.  Let $\lll$ be a set of representatives of
  these equivalence classes, restricted to those that are composed of edges with lengths in $\Q(\beta)$, so that they are commensurate with the edges of $\alpha$.  Without loss of generality we choose
  these representatives to be on the positive real axis with an
  endpoint at the origin.  For each $L \in \lll$, the endpoints of the
  lower tiles of $L$ create a partition $0 = a_0 < a_1 < ... <
  a_{k_\ell} = |L|$, where $k_\ell$ is the number of lower tiles of
  $L$ and $|L|$ is the length of $L$.  Any such $a_i$ for any initial
  segment $L$ will be called a {\em prefix} of $L$.  Similarly, the
  endpoints of the upper tiles of $L$ create a partition $0=b_0<b_1<
  ... < b_{k_u}=|L|$ where $k_u$ is the number of upper tiles of $L$.
  A {\em suffix} of an initial segment $L$ is $|L| - b_i$, where $b_i$
  is in the partition of $L$.  In the same manner, we may compute
  prefixes of $S_\ell(e)$ for $e \in \eee_\ell$ and suffixes of
  $S_u(e)$ for $e \in \eee_u$,  restricting again to edges with lengths in $\Q(\beta)$.  Since $|S_\ell(e)| = |S_u(e)| = \beta
  |e|$, the prefixes will be partition elements of $[0,\beta |e|]$,
  and the suffixes will be $\beta|e|$ minus partition elements of $
  [0, \beta |e|]$.  We collect all such prefixes into a prefix
  set $\qqq_{px}$ and such suffixes into a suffix set
  $\qqq_{sx}$. It is important to note that both $\qqq_{px}$ and
  $\qqq_{sx}$ are finite.
  
  Since each $P\in \ppp$ is fractagonal, there are only a finite
  number of straight edges. Thus there is a $q \in\N$ so that all
  straight edge lengths in the class of the simple adjacency $\alpha$
  are in $\frac{1}{q}\Z[\beta]$.  This implies
  that the lengths of all initial segments, and the lengths of all
  prefixes and suffixes in $\qqq_{px}$ and $\qqq_{sx}$ are also in
  $\frac{1}{q}\Z[\beta]$. Since offsets are the difference between a
  sum of upper edges and a sum of lower edges,
  $\ooo\subset\frac{1}{q}\Z[\beta]$.
    
  We have already seen that $\ooo$ lies in a segment $L'$ that is the
  image $\lambda^n L''$ for some initial segment $L''$.  Let $L_0$ be
  the representative in $\lll$ of $L''$.  The tiles in $\alpha$, which
  intersect $L' = \lambda^n\cdot L''$, are equivalent to lower and
  upper tiles $T_\ell$ and $T_u$ intersecting $\beta^n\cdot L_0$.  In
  fact the edges of $T_\ell$ and $T_u$ are lower and upper edges that
  are edges in $S_\ell(L_0)$ and $S_u(L_0)$ when we consider $L_0$ in
  terms of its lower and upper tiles.  The endpoints of edges in
  $\alpha$ used to compute $\ooo$ correspond to points $a \in T_\ell$
  and $b \in T_u$ (see Figure \ref{upper.lower.offset}), so to compute
  $\ooo$ we need to compute $|b - a|$.  Without loss of generality
  assume $a \le b$.
  
  \begin{figure}[h]
\epsfig{figure=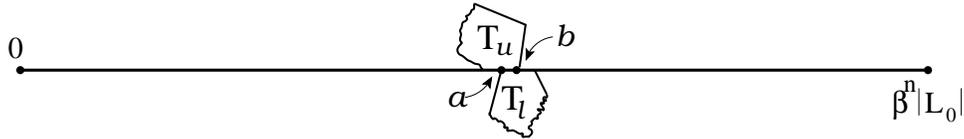,width=5in}  
\caption{The simple adjacency $\alpha$ is equivalent to $T_u \cup T_{\ell}$.}
  \label{upper.lower.offset}
  \end{figure}
  
  We may compute the left endpoint $a$ in terms of powers of $\beta$
  applied to elements of $\qqq_{px}$; we do the same for the upper
  endpoint $b$ using powers of $\beta$ applied to elements of
  $\qqq_{sx}$.  (See Figure \ref{sum.powers}).  Let $p_n \in
  \qqq_{px}$ be the largest prefix of $L_0$ for which $\beta^n p_n \le
  a$.  If equality holds, we have computed $a$, otherwise we continue.
  We denote by $e_{n-1}$ the lower edge immediately following $p_n$ in
  $L_0$.  This edge will, after application of $S_l^{n}$, contain $a$
  (otherwise $p_n$ wasn't maximal).  So there is a largest prefix of
  $e_{n-1}$ which we call $p_{n-1}\in \qqq_{px}$ for which $\beta^n
  p_n + \beta^{n-1} p_{n-1}\le a$.  If equality holds, we are done and
  the remainder of the prefixes will be the zero prefix.  If not, we
  have that $p_{n-1}$ is a prefix of $S_\ell(e_{n-1})$ for some edge
  $e_{n-1} \in \eee_\ell$, so the lower edge immediately following
  $p_{n-1}$ in $S_\ell(e_{n-1})$ will, after application of $S^{n-1}$,
  contain $a$.  We can continue in this fashion, choosing prefixes
  $p_{n-2}, ..., p_0$ until we obtain that $a = \sum_{k = 0}^{n}
  \beta^k p_k$.
  
  Similarly we can compute the value of $b$.  Let $s_n\in \qqq_{sx}$
  be the largest suffix of $L_0$ for which $b \le \beta^n(|L_0|-
  s_n)$.  If equality holds we are done and all subsequent suffixes
  are the zero suffix.  If not we select the edge $e_{n-1}$
  immediately preceding $|L_0|-s_n$ and find its suffix $s_{n-1}$ so
  that $b \le \beta^n(|L_0|- s_n)-\beta^{n-1}s_{n-1}$.  We continue
  choosing suffixes $s_{n-2}, ..., s_0$ until we obtain that $b =
  \beta^n|L_0| -\sum_{k = 0}^{n } \beta^k s_k.$

\begin{figure}[ht]
\epsfig{figure=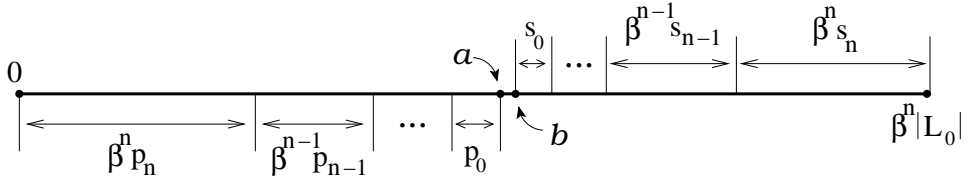,width=5in}
\caption{The computation of $a$ and $b$ as seen in $[0,\beta^{n} |L_0|].$}
\label{sum.powers}
\end{figure}

Thus the offset distance $\ooo$ is given by 
\begin{equation}\label{length.oi}
\ooo =  \beta^n|L_0| + \sum_{k = 0}^{n} \beta^k(-s_k - p_k).
\end{equation}
Our interest, then, is in all possible combinations of suffixes
and prefixes. We collect these into a {\em difference set}
$\ddd=\{-s-p \text{ such that } s \in \qqq_{sx}, p \in \qqq_{px}\}$ which we must show  is finite. 

If $\beta$ is an integer, then 
the result follows immediately since $\ooo$ is bounded by the length
of the longest straight edge and is an element of $\frac{1}{q} \Z$.

If $\beta$ is not an integer, we need to use the ideas and notation of
Section \ref{OneDResults}.  Recall that the map $\phi: \C^d \to \C$ is
a bijection between $\Q^d$ and $\Q(\beta)$, and that multiplication by
the matrix $A$ in $\Q^d$ is equivalent to multiplication by $\beta$ in
$\Q(\beta)$.  We can find preimages $\vecw_k \in \Q^d$ such that
$\phi(\vecw_k) = -s_k - p_k$ for $k = 0, ..., n $.  Since elements of
$\ddd$ are in $\frac{1}{q}\Z[\beta]$, we have that each $\vecw_k $ must
be in $\frac{1}{q}\Z^d$.  Equation (\ref{length.oi}) becomes
\begin{equation}
\ooo = \beta^n|L_0| +\sum_{k = 0}^{n} \beta^k \phi(\vecw_k)
= \beta^n|L_0| +\phi\left( \sum_{k = 0}^{n}A^k \vecw_k \right).
\end{equation}

We write the preimages of $\ddd$ restricted to $\Q^d$ in terms of the
eigenvectors $\vecv_1, ..., \vecv_d$.  Since $\phi$ projects onto
$\vecv_1$ sending all other eigenvectors to 0, and since $\ooo$ is
necessarily less than the length of the longest edge, there is an $m_1
> 0$ that is an upper bound on the modulus of the component of $\phi^{-1}(\ooo) \cap \Q^d$
in the $\vecv_1$ direction.  For any $j = 2, ..., d$, denote by $m_j$
the maximum modulus of the component in the $\vecv_j$ direction of any
preimage of the difference set $\ddd$ restricted to $\Q^d$.  The modulus of the component of
$\sum_{k = 0}^{n}A^k \vecw_k$ in the $\vecv_j$ direction is less than
or equal to
$$\sum_{k = 0}^{n} |\beta_j^k| m_j < \sum_{k = 0}^{\infty} |\beta_j^k| m_j
= \frac{m_j}{1 - |\beta_j|},$$
since $|\beta_j| < 1$ by assumption when
$j = 2, ..., d$.  So we see that the vectors which project onto $\ooo$
must lie in a bounded region.  There are only a finite number of
elements of $\frac{1}{q}\Z^d$ in a bounded region, so there are only a
finite number of offsets $\ooo$, and hence only a finite number of
simple adjacencies $A$.  The tiling $x$ is locally finite.
\end{proof}

\noindent{\bf Comment.\ }
A complex algebraic integer $\lambda$ is called a {\em complex Pisot
  number} if $|\lambda|>1$ and all algebraic conjugates
$\lambda'\not=\overline \lambda$ satisfy $|\lambda'|<1$ (see
\cite{Sol.self.similar}, \cite{Robinson.ams}). If 
$\lambda$ is complex Pisot, $\beta=|\lambda|$ may not be real Pisot
(see \cite{Robinson.ams} for an example). We do not know whether
Theorem~\ref{thm.flc}  always holds in the case of a complex Pisot
  linear  expansion $\lambda$.

\subsection{The Figure~\ref{np.subs} example is not locally finite}
\label{Intro.example.non.flc}

Let $a$ denote the large square prototile from Figure \ref{np.subs}.
Since the substitution of any tile includes a tile of type $a$, it
follows that a block of type $S^n(a)$ must appear in $S^{n+1}(e)$ for
any tile type and for any $n = 1, 2, ...$ and so every tiling $x \in
X$ must contain blocks of type $S^n(a)$ for all $n = 1, 2, ...$.

The circled point in Figure \ref{iterate.misfit} appears inside the
block $S(a)$.  This point represents a vertex of the tiles above it
and it intersects the upper boundary of the tile below it at a point
$\beta - 1$ units from the left endpoint of that boundary.  Each time
we apply the substitution, the circled point determines the offset
between the tile below it and the two tiles above it.  This offset is
completely determined by $F$.  If we can show that the orbit of $\beta
- 1$ is infinite then this implies that there are an infinite number
of offsets and hence the tiling is not locally finite.  It is
convenient to compute this orbit for $f$ rather than $F$.

\begin{figure}[ht]
\epsfig{figure=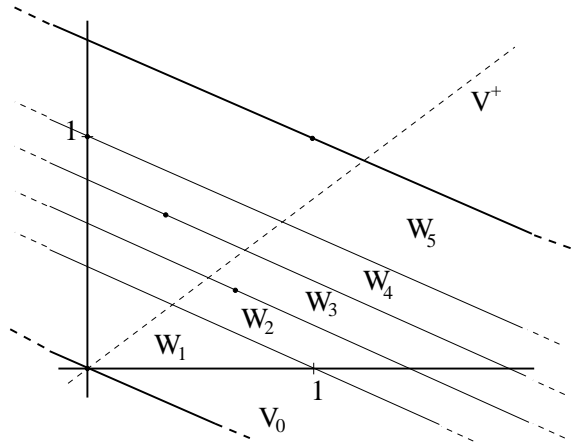,width=3in}
\caption{The regions $W_j$ in the domain $W_E$ of $f$.}
\label{fff.ex}
\end{figure}

Example \ref{BetaEx} gives the algebraic generalized
$\beta$-transformation associated with the substitution on upper edges
of tiles.  Noting that the eigenvalues are real, we compute the
corresponding action $f$ in $W_E \subset \R^2$.  The companion matrix
of $\beta$ is
$$A= \left[ \begin{array}{cc} 0&1\\1&3 \end{array}\right].$$
We choose as the generator of $V^+$ the vector $\vecv = (-\beta_2, 1)$ and as
generator of $V_0$ we take $\vecv_2 = (-\beta, 1)$.  The vectors
mapping onto $y_1, y_2, y_3, y_4$ and $y_5$ are given by $\vecw_1 =
(0,0), \vecw_2 = (0,0), \vecw_3 = (1,0), \vecw_4 = (2,0),$ and $
\vecw_5 = (3,1)$.  The sets $W_j$ have
boundaries given by $V_0 + \phi^{-1}(x_{j - 1})$ and $V_0 +
\phi^{-1}(x_j)$, where the inverse is taken in $\Q^2$.  It is not
difficult to check that the appropriate translation vectors are
$(0,0), (1,0), (2/3, 1/3), (1/3, 2/3), (0,1)$ and $(1,1)$.  (In this
example one can also find these by looking at $A^{-1} \vecw_j$.)  The
regions $W_j$ are as shown in Figure \ref{fff.ex}, along with the main
eigenspace $V^+$.

We now establish an $f$-orbit which enters the escape region and thus
tends to infinity.  By our theorem, this orbit corresponds to an
infinite $F$-orbit.  We begin by computing the escape region.  Looking
in the proof of Theorem \ref{orbit.type} (\ref{nonpisot.infinite}),
we see that we must choose $r_2$ so that $$
r_2 > \frac{|w_{j,
    2}|}{|\beta_2| - 1} \quad \text{ for all } \quad j = 1, 2, 3, 4.$$
We recall that $\vecw_1 = (0,0), \vecw_2 = (1,0), \vecw_3 = (2,0),$
and $ \vecw_4 = (3,1)$ and compute their components in the $V_0$
direction:
$$w_{1,2} = 0, \quad w_{2,2}=\frac{-1}{\sqrt{13}}, \quad
w_{3,2}=\frac{-2}{\sqrt{13}}, \quad w_{4,2} = \frac{-7 +
  \sqrt{13}}{2\sqrt{13}}.$$
The largest of these in modulus is
$|w_{3,2}| \approx .5547$.  We see that $$
r_2 > \frac{|w_{3,
    2}|}{|\beta_2| - 1} \approx 1.8321, $$
so we simply choose $r_2 =
2$ to define $\resc$.

We will show that $\vecw =\phi^{-1}(\beta - 1) = (-1, 1)$ enters the
escape region by computing the first several elements of its orbit.
(This was made easier by using a cobwebbing program in Matlab to
determine in which $W_j$ the point was at each step.)  Now $(-1,1) \in
W_1$, so $f(-1,1) = A(-1,1) - \vecw_1 = (3,0).$ One can compute that
$(3,0) \in W_4$, so $f(3,0) = A(3,0) - \vecw_4 = (-3,2)$.  Continuing
in this fashion the next several iterates are $(5, -1), (-6, 3), (9,
-3), (-11,6), (15, -6)$.  The orbit enters the escape region at
$(-6,3)$ and the iterates begin to grow quite quickly.  It is
instructive and surprising to draw the region $W_E$, which is infinitely
long and quite thin, and plot these points which must remain inside
but which bounce further and further away from the origin.

\bibliographystyle{amsplain}

\end{document}